\documentclass[leqno]{article} 
\usepackage{latexsym}
\input epsf.sty
\usepackage{amssymb}
\begin{document}
\newtheorem{proposition}{Proposition}[section]
\newtheorem{application}{Application}[section]
\newcommand{\Var}{\mathop{\rm Var}}
\newcommand{\ds}{\displaystyle}
\newcommand{\Prob}{\mathop{\rm P}}
\newcommand{\barF}{\overline{F}}
\newcommand{\EE}{\mathop{\rm E}}
\newcommand{\Stati}{\hbox{$\cal S$}}
\newcommand{\Reali}{\mathbb{R}}
\newcommand{\fine}{\hfill $\Box$}             
\newcommand{\qed}{\hfill\rule{2mm}{2mm}}      
\def\eq#1{~{\rm (\ref{equation:#1})}}
\def\LEQ#1{\mathop{{\leq}_{\hbox{\scriptsize\rm #1}}}}
\newenvironment{proof}{\begin{trivlist}
\item[\hspace{\labelsep}{\bf\noindent Proof. }]
}{\fine\end{trivlist}}
\title{ON THE EFFECT OF RANDOM ALTERNATING PERTURBATIONS 
ON HAZARD RATES\footnote{Paper appeared in Sci.\ Math.\ Jpn.\ 64 (2006), no. 2, 381--394.}}
\author{\sc Antonio Di Crescenzo {\rm and} Barbara Martinucci \\
\rm Dipartimento di Matematica e Informatica, Universit\`a di Salerno \\
\rm Via Ponte don Melillo, I-84084 Fisciano (SA), Italy.}

\date{\empty}

\maketitle
\begin{abstract}
We consider a  model for systems perturbed by dichotomous noise, in which 
the hazard rate function of a random lifetime is subject to additive 
time-alternating perturbations described by the telegraph process. This leads 
us to define a real-valued continuous-time stochastic process of alternating type 
expressed in terms of the integrated telegraph process for which we obtain the 
probability distribution, mean and variance. An application to survival analysis 
and reliability data sets based on confidence bands for estimated hazard rate 
functions is also provided. 



\end{abstract}
%
\section{Introduction}\label{section:1}
Large attention has been given recently in the physical and mathematical literature to 
stochastic systems perturbed by dichotomous noise, such as those related to Brownian 
motion or Brownian motors (see, for instance, Bena {\em et al.}\ \cite{BBKK02}, \cite{BBKK03} 
and Porr\`a {\em et al.}\ \cite{PoRoMa96}). The relevance of Markovian dichotomous noise in 
biological systems has been also pointed out. See, for instance, Laing and Longtin \cite{LaLo01} 
where the benefic role of noise is investigated within the context of mathematical neuroscience. 
In this paper we aim to study the effect of dichotomous noise to the hazard rate function 
of random lifetimes.
\par
It is well known that in many fields related to survival analysis and reliability theory 
the study of hazard rates plays a very important role (see the classical book of 
Barlow and Proschan \cite{BaPr96}, for instance). Specifically, here we are interested 
to hazard rate functions that are realizations of stochastic processes. 
Such types of hazard rates arise in some models that were recently proposed to 
describe doubly random phenomena such as lifetimes of devices operating in 
random environments. See, for instance, Kebir \cite{Ke91} and the included 
references for a model in which the hazard rate is a functional of a stochastic 
process that describes the random variability of the environment. 
A similar model is treated in Di Crescenzo and Pellerey \cite{DiCrPe98}, 
where the hazard function is the realization of a non-decreasing 
stochastic process with independent increments. Other models characterized by 
stochastic hazard rate functions were also discussed by 
\={E}\v{z}ov and Jadrenko \cite{EzJa67} and by Yadrenko and Zhegri\u\i \ \cite{YaZh95}. 
Furthermore, a stochastic model in which 
a constant hazard rate (describing the fractional dissolution rate of a drug potion) 
is corrupted by a white noise was recently analyzed by L\'ansk\'y and L\'ansk\'a 
\cite{LaLa00} and by L\'ansk\'y and Weiss \cite{LaWe01} within the theory of drug 
dissolution. Other applications of stochastic hazard rates were proposed 
in the actuarial literature and in studies of population biology and human aging 
(see Milevsky and Promislow \cite{MiPr01} and references therein). 
All these stochastic models are suitable to describe phenomena characterized 
by an intrinsic randomness, which is evidenced at the hazard rate level. 
\par
Along the line traced by the above mentioned papers, herewith we propose a new 
model to describe systems characterized by failure rates which are subject 
to random perturbations expressed as dichotomous noise. We aim to discuss 
a model that incorporates a more realistic kind of noise characterized by a constant  
intensity, and that at the same time presents the feature of noise with alternating 
behavior that is frequently found in biological systems. For instance, see the papers by 
Buonocore {\em et al.}\ \cite{BuDiCrDiNa2002} and Di Crescenzo {\em et al.}\ \cite{DiCrMaPi2005}, 
dealing with neuronal dynamics subject to alternating inputs. Stochastic models characterized 
by alternating behaviour are also used to describe the cell motility on one-dimensional 
and two-dimensional state-spaces (see, for instance, Lutscher \cite{Lu03}, 
Stadje \cite{St87}, \cite{St88}, and references therein). 
\par
Suppose that (non-negative) lifetime $T$ has distribution function
\begin{equation}
   F(t)=1-\exp\left\{-\int_{0}^{t}r(s)\,{\rm d}s\right\},
   \qquad t\geq 0,
   \label{equation:20}
\end{equation}
so that its hazard rate function is 
\begin{equation}
   r(t)=\lim_{\Delta t \rightarrow 0^+} \frac{\Prob(T\leq t+\Delta t\,|\,T>t)}{\Delta t},
   \qquad t\geq 0.
   \label{equation:6}
\end{equation}
We assume that $r(t)>c$, where $c$ is a positive constant, and that $r(t)$ is subject 
to additive time-alternating perturbations described by the well-known telegraph  
process (see, for instance, Orsingher \cite{Or90} or Beghin {\em et al.}\ \cite{BeNiOr01}). 
This assumption consists of applying the substitution 
$$
 r(t)\;\rightarrow \;r(t)+V(0) (-1)^{N(t)},\qquad t\geq 0 
$$
in the right-hand-side of\eq{20}, where $V(0)\in\{-c,c\}$ and $\{N(t); t\geq 0\}$ is 
a Poisson process independent of $T$. We are thus led to defining a real-valued 
continuous-time stochastic process $\{X(t); t\geq 0\}$, where $X(t)$ is expressed 
in terms of the integrated telegraph process:
\begin{equation}
   X(t)=1-\exp\left\{-\int_{0}^{t}[r(s)+V(0) (-1)^{N(s)}] \,{\rm d}s\right\},
   \qquad t\geq 0.
   \label{equation:19}
\end{equation}
\par
Here, we present some basic results on process\eq{19}. First of all, 
in Section \ref{section:2} we obtain the moment generating function of the 
integrated telegraph process, which is involved in the expression of 
mean and variance of $X(t)$. Such moments and the probability distribution 
of $X(t)$ are given in Section \ref{section:3}. This distribution 
has a discrete as well as an absolutely continuous component. 
Finally, in Section \ref{section:4} we use an asymptotic confidence 
band for estimated hazard rates in two case-studies where our model 
provides adequate fit to two data sets taken from the survival analysis 
and reliability theory literature. 
\par
We point out that the results obtained in this paper can be extended to the case where 
there are two different values $V(0)\in\{-c_2,c_1\}$ and $N(t)$ is an alternating Poisson 
process characterized by two rates $\lambda_1$ and $\lambda_2$. Finally, we point 
out that stochastic processes obtained by transformations of the integrated 
telegraph process are of interest in stochastic modelling literature. 
For example, Di Crescenzo and Pellerey \cite{DiCrPe02} use a geometric 
telegraph process to model the price dynamics of risky assets.
%
\section{Moment generating function of integrated telegraph process}\label{section:2}
%
Let $\{W(t); t\geq 0\}$ denote the well-known integrated telegraph process,
defined by
\begin{equation}
    W(t)=V(0)\,\int_0^t (-1)^{N(s)}\,{\rm d}s,
    \qquad t\geq 0,
    \label{equation:1}
\end{equation}
where $\{N(t); t\geq 0\}$ is a homogeneous Poisson process with intensity $\lambda$ 
and $V(0)$ is a random variable independent from $N(t)$ such that
\begin{equation}
\Prob\{V(0)=\pm c\}=\frac{1}{2},
    \label{equation:12}
\end{equation}
with $c>0$. 
The distribution of\eq{1} is characterized by discrete components concentrated
at $\pm c\,t$ and by an absolutely component on $(-c\,t,c\,t)$. Orsingher~\cite{Or90}, 
shows that the discrete component satisfies 
\begin{equation}
    \Prob\{W(t)=c\,t\}=\Prob\{W(t)=-c\,t\}=\frac{1}{2}\,{\rm e}^{-\lambda\,t},
    \qquad t\geq 0,
    \label{equation:2}
\end{equation}
and that, for $t\geq 0$, $-ct<x<ct$, the continuous component has density 
\begin{equation}
    {\partial\over \partial x}\Prob\{W(t)\leq x\}
    =\frac{1}{2c}\,{\rm e}^{-\lambda\,t}
    \left[\lambda\,I_0\left(\frac{\lambda}{c}\sqrt{c^2t^2-x^2}\right)
    +{\partial\over \partial t}I_0\left(\frac{\lambda}{c}\sqrt{c^2t^2-x^2}\right)\right],
    \label{equation:3}
\end{equation}
where
\begin{equation}
 I_0(x)=\sum_{k=0}^{+\infty}\frac{(x/2)^{2k}}{(k!)^2}
    \label{equation:24}
\end{equation}
denotes the modified Bessel function of order $0$. 
\begin{proposition}
For all $s\in\Reali$ and $t\geq 0$ the moment generating function of $W(t)$ is
\begin{eqnarray}
    && M(s,t):=E\left[{\rm e}^{s\,W(t)}\right]
    \label{equation:4} \\
    && \hspace{1.4cm}
    ={\rm e}^{-\lambda\,t}
    \left[\cosh\left(t\sqrt{\lambda^2+s^2c^2}\right)
    +{\lambda\over \sqrt{\lambda^2+s^2c^2}}\,
    \sinh\left(t\sqrt{\lambda^2+s^2c^2}\right)\right].
    \nonumber
\end{eqnarray}
\end{proposition}
\begin{proof}
   From Eqs.\eq{2} and\eq{3} it follows
\begin{eqnarray}
    M(s,t) \!\!\!\!
    &=& \!\!\!\! {{\rm e}^{-\lambda t}\over 2}\left({\rm e}^{sct}+{\rm e}^{-sct}\right) 
    \label{equation:5} \\
    & & \!\!\!\!  
    +{{\rm e}^{-\lambda t}\over 2c}\int_{-ct}^{ct}{\rm e}^{sx}
    \left[\lambda\,I_0\left({\lambda\over c}\sqrt{c^2t^2-x^2}\right)
    +{\partial\over\partial t}I_0\left({\lambda\over c}\sqrt{c^2t^2-x^2}\right)\right]{\rm d}x
    \nonumber \\
    &=& \!\!\!\! {{\rm e}^{-\lambda t}\over 2c}\left[\lambda\,Q(s,t)
    +{\partial\over\partial t}Q(s,t)\right],
    \qquad s\in\Reali,\;t\geq 0,
    \nonumber 
\end{eqnarray}
where we have set
$$
    Q(s,t):=\int_{-ct}^{ct}{\rm e}^{sx}\,
    I_0\left({\lambda\over c}\sqrt{c^2t^2-x^2}\right){\rm d}x.
$$
Making use of Eq.\ (25), from Orsingher~\cite{Or90} we obtain:
$$
    \int_{-ct}^{ct}{\rm e}^{sx}{\partial^2\over\partial t^2}
    I_0\left({\lambda\over c}\sqrt{c^2t^2-x^2}\right){\rm d}x
    =\int_{-ct}^{ct}{\rm e}^{sx}c^2 {\partial^2\over\partial x^2}
    I_0\left({\lambda\over c}\sqrt{c^2t^2-x^2}\right){\rm d}x
    +\lambda^2\,Q(s,t).
$$
A two-fold integration by parts shows that 
$$
    {{\rm d}^2\over {\rm d}t^2}Q(s,t)=(\lambda^2+s^2c^2)\,Q(s,t).
$$
Solving this equation with initial conditions $Q(s,0)=0$ and
$\ds{{\rm d}\over {\rm d}t}Q(s,t)\big|_{t=0}=2c$ we have:
$$
    Q(s,t)={c\over \sqrt{\lambda^2+s^2c^2}}\left[{\rm e}^{t\sqrt{\lambda^2+s^2c^2}}
    -{\rm e}^{-t\sqrt{\lambda^2+s^2c^2}}\right],
    \qquad s\in\Reali,\;t\geq 0.
$$
Using this formula in the right-hand-side of\eq{5}, expression\eq{4} 
finally follows.
\end{proof}
\par
It should be noticed that\eq{4} could also be obtained from the initial-value
problem for the telegraph equation
$$
    \left\{
    \begin{array}{l}
    \ds{\partial^2\over\partial t^2}p+2\lambda\ds{\partial\over\partial t}p
    =c^2\ds{\partial^2\over\partial x^2}p \\
    \vspace{-0.2cm}\\
    p(x,0)=\delta(x) \\
    \vspace{-0.2cm}\\
    \ds{\partial\over\partial t}p(x,t)\Big|_{t=0}=0,
    \end{array}
    \right.
$$
where $\delta(x)$ is the Dirac delta function. Indeed, $M$ is solution of
$$
    \left\{
    \begin{array}{l}
    \ds{{\rm d}^2\over{\rm d}t^2}M+2\lambda\ds{{\rm d}\over{\rm d}t}M
    =s^2c^2M \\
    \vspace{-0.2cm}\\
    M(s,0)=1 \\
    \vspace{-0.2cm}\\
    \ds{\partial\over\partial t}M(s,t)\Big|_{t=0}=0.
    \end{array}
    \right.
$$
\section{Results}\label{section:3}
Let $T$ be a non-negative random variable, with absolutely continuous
distribution function\eq{20}, hazard rate function\eq{6}, and support
$$
    {\cal S}:=[0,\ell),
    \qquad \hbox{where } \;\;\ell=\sup\{t\geq 0: F(t)<1\}\in(0,+\infty].
$$
Throughout the paper we will assume that  $r(t)>c$ for all $t\in{\cal S}$,
where $c$ is a positive constant. Hence,  
$$
    F(0)=0, \qquad F(t)>1-{\rm e}^{-c t} \quad \hbox{for all }t>0,
$$
so that 
$$
    T\LEQ{st}{\rm Exp}(c),
$$
where `$\LEQ{st}$' denotes the usual stochastic order (see, for instance,
Shaked and Shantikumar~\cite{ShSh94}) and ${\rm Exp}(c)$ denotes an
exponentially distributed random variable with mean $c^{-1}$. 
\par
Let us now consider the stochastic process $\{X(t); t\in{\cal S}\}$, with
$X(0)=0$ and
\begin{equation}
    X(t)=1-\barF(t)\,{\rm e}^{-W(t)},
    \qquad t\in{\cal S},
    \label{equation:8}
\end{equation}
where $W(t)$ has been defined in\eq{1} and where $\barF(t)=1-F(t)$ denotes the survival 
function of $T$. Note that process\eq{8} is equivalent to process\eq{19}. 
Hence, for any fixed $t$, $X(t)$ describes the probability that a random lifetime 
is not larger than $t$. The system hazard rate is random and is given by 
$r(t)+V(0) (-1)^{N(t)}$, where the distribution of $V(0)$ is specified in\eq{12}. 
As pointed out in Section 1, this hazard rate accounts for random perturbations on $r(t)$ 
that occur according to a telegraph process. The sample-paths of $X(t)$ are 
absolutely continuous distribution functions that approach\eq{20} as $c$ goes to $0$. 
As an example, Figure 1 shows (a) two sample-paths of $W(t)$, and (b) the 
corresponding sample-paths of $X(t)$, characterized by distribution function 
%
\begin{figure}[t] 
\begin{center}
\epsfxsize=8.5cm
\centerline{\epsfbox{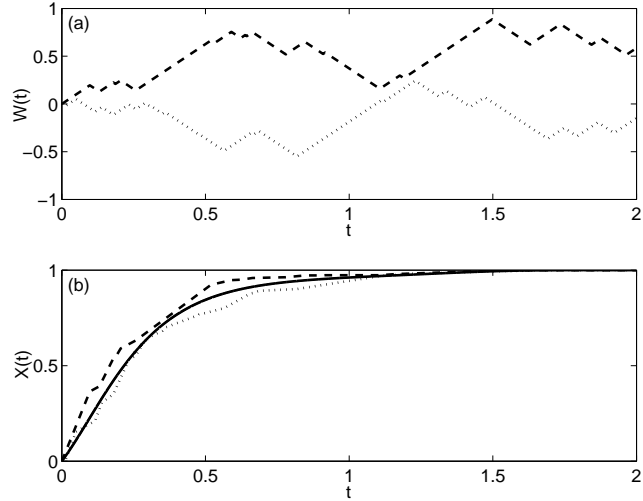}}
\caption{\small (a) Two sample-paths of process\eq{1} for $c=2$ and $\lambda=15$; 
(b) the corresponding sample-paths of\eq{8}, with $F(t)$ and $r(t)$ 
defined in\eq{22} and\eq{23}, respectively, with $\alpha=15$ and $\beta=0.001$; 
the solid line shows $F(t)$.}
\label{fig:1}
\end{center}
\end{figure}
%
\begin{equation}
 F(t)=1-\exp\left\{-\left[{\alpha\over 4}\,t^4-{2\alpha\over 3}\,t^3+{\alpha\over 2}\,t^2+(c+\beta)\,t\right]\right\}, 
 \qquad t\geq 0,
    \label{equation:22}
\end{equation}
where 
\begin{equation}
 r(t)=\alpha t(t-1)^2+c+\beta, \qquad t\geq 0, \quad (\alpha,\beta> 0).
    \label{equation:23}
\end{equation}
\par
Next, we derive the probability distribution of $X(t)$. Since 
$\Prob\{-ct\leq W(t)\leq ct\}=1$, $t\geq 0$, from Eq.\eq{8} we have:
$$
 \Prob\{a(t)\leq X(t)\leq b(t)\}=1, \qquad t\in{\cal S},
$$
and   
$$
 a(t)<F(t)<b(t)\quad \hbox{for all }t\in(0,\ell),
$$
where
\begin{equation}
    a(t)=1-\barF(t)\,{\rm e}^{ct},  \qquad
    b(t)=1-\barF(t)\,{\rm e}^{-ct}, \qquad t\in{\cal S}.
    \label{equation:9}
\end{equation}
From\eq{9}, we note that $a(t)$ and $b(t)$ are distribution functions 
with hazard functions $r(t)-c$ and $r(t)+c$, respectively. 
Furthermore, there results:
$$
  \lim_{t\to \ell}a(t)=1-{\rm e}^{-\nu},
  \qquad
  \lim_{t\to \ell}b(t)=1,
$$
where  
$$
\nu=\int_0^{\ell}[r(s)-c]\,{\rm d}s.
$$
We note that $\nu=+\infty$ when $\ell<+\infty$, whereas $\nu\leq +\infty$ if $\ell=+\infty$. 
Moreover, by setting $D(t):=b(t)-a(t)$, from\eq{9} we obtain:
$$
  D(t)=2\,\barF(t)\,\sinh(ct), \qquad t\in{\cal S},
$$
with $D(0)=0$ and $\ds\lim_{t\to\ell}D(t)={\rm e}^{-\nu}$, with $D(t)$   
non-decreasing in $t$ if
$$
  r(t)\leq c\,\coth(ct).
$$
For instance, $D(t)$ is shown in Figure 2 for three different choices of $r(t)$. 
The first choice shows a case in which $D(t)$ is bimodal. In the third choice we 
have $\ds\lim_{t\to +\infty}D(t)=0.268$, whereas such limit vanishes in the first 
two choices. 
%
\begin{figure}[t] 
\begin{center}
\epsfxsize=8.5cm
\centerline{\epsfbox{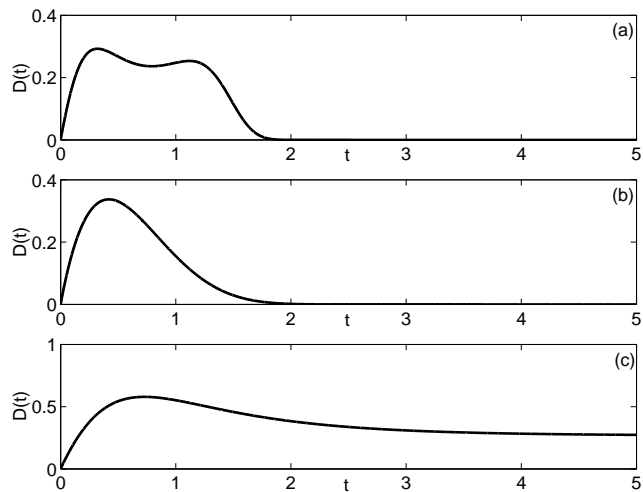}}
\caption{\small Plots of $D(t)$ in 3 different cases: 
(a) $r(t)$ given by Eq.\eq{23}, with $c=1$, $\lambda=\alpha=15$ and $\beta=0.001$; 
(b) $r(t)=1+{\rm e}^t$, $t\geq 0$, with $c=\lambda=1$; 
(c) $r(t)=1+3(1-{\rm e}^{-t})/({\rm e}^t+{\rm e}^{-t})$, $t\geq 0$, with $c=\lambda=1$.}
\label{fig:2}
\end{center}
\end{figure}
%
\par
We shall now determine the probability distribution of $X(t)$. Like as for the integrated 
telegraph process, the distribution of $X(t)$ consists of a discrete component on 
$a(t)$ and $b(t)$ and an absolutely continuous component inside of $\big(a(t),b(t)\big)$. 
\begin{proposition}
For all $t\in{\cal S}$ the discrete component of the distribution of $X(t)$ is: 
\begin{equation}
    \Prob\{X(t)=a(t)\}=\Prob\{X(t)=b(t)\}=\frac{1}{2}\,{\rm e}^{-\lambda\,t}.
    \label{equation:21}
\end{equation}
Moreover, for all $x\in\big(a(t),b(t)\big)$ and $t\in{\cal S}$ the continuous 
component has density 
\begin{equation}
    f(x,t)
    =\frac{1}{2c}\,\frac{1}{1-x}\,{\rm e}^{-\lambda\,t}
    \left[\lambda\,I_0\left(\frac{\lambda}{c}\sqrt{u(x,t)}\right)
    +{\partial\over \partial t}I_0\left(\frac{\lambda}{c}\sqrt{u(x,t)}\right)\right],
    \label{equation:17}
\end{equation}
where
\begin{eqnarray}
    && u(x,t):=c^2t^2-\ln^2\frac{\barF(t)}{1-x}
    \label{equation:13} \\
    && \hspace{1.25cm}
    =\ln\left(\frac{1-a(t)}{1-x}\right)\,\ln\left(\frac{1-x}{1-b(t)}\right),
    \nonumber
\end{eqnarray}
with $a(t)$ and $b(t)$ defined in\eq{9}.
\end{proposition}
\begin{proof}
Due to\eq{8}, the distribution function of $X(t)$ can be expressed as
$$
 \Prob\{X(t)\leq x\}=\Prob\left\{W(t)\leq \ln{\barF(t)\over 1-x}\right\},
 \qquad t\in{\cal S},\;\; a(t)\leq x\leq b(t).
$$
Hence, recalling\eq{2} and\eq{3}, Eqs.\eq{21} and\eq{17} follow.
\end{proof}
\par
The expression of density\eq{17} is similar to that of the underlying integrated 
telegraph process. Unlike the latter process, however, the support of $f(x,t)$ 
is bounded when $t$ increases. Moreover, the expression for $f(x,t)$ reflects the 
form of interval $(a(t),b(t))$. Indeed, from\eq{13} it is not hard to note that 
$u(x,t)=0$ when $x=a(t)$ and $x=b(t)$, whereas $u(x,t)>0$ when $x\in (a(t),b(t))$. 
Furthermore, recalling\eq{9} and\eq{24}, for $t\in\Stati$ we have: 
$$
 \lim_{x\downarrow a(t)}f(x,t)
 ={\lambda\,{\rm e}^{-(\lambda+c)t}\over 2c\,\barF(t)}\left(1+{\lambda t\over 2}\right),
 \qquad 
 \lim_{x\uparrow b(t)}f(x,t)
 ={\lambda\,{\rm e}^{-(\lambda-c)t}\over 2c\,\barF(t)}\left(1+{\lambda t\over 2}\right).
$$
Figure 3 shows for instance $f(x,t)$ for three different choices of $t$. 
%
\begin{figure}[t] 
\begin{center}
\epsfxsize=8.5cm
\centerline{\epsfbox{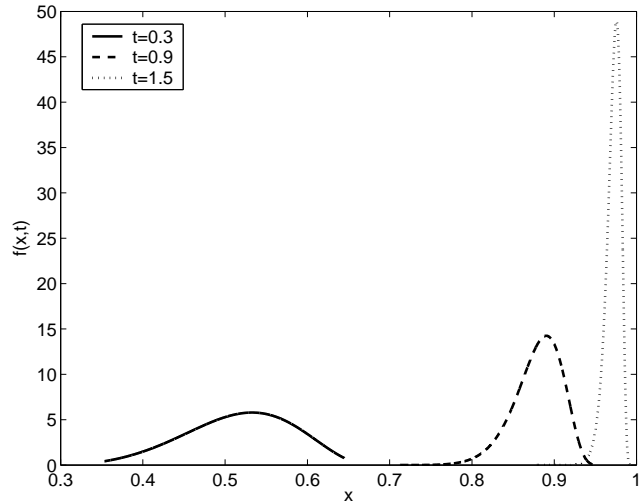}}
\caption{\small Plot of density\eq{17} for three choices of $t$ 
and for $x$ ranging within $(0,1)$; 
$r(t)$ and all parameters are chosen as in case (a) of Figure 2, i.e.\ $c=1$, 
$\lambda=\alpha=15$, $\beta=0.001$ and $r(t)$ is given by Eq.\eq{23}.}
\label{fig:3}
\end{center}
\vspace{2cm}
\end{figure}
%
\begin{proposition}
For all $t\geq 0$ mean and variance of $X(t)$ are:
\begin{eqnarray}
&& \hspace{-0.5cm}
    \EE\left[X(t)\right]=1-\barF(t)\,M(-1,t),
    \label{equation:10} \\
&& \hspace{-0.5cm}
    \Var\left[X(t)\right]=\barF^2(t)\,\left\{M(-2,t)-[M(-1,t)]^2\right\},
    \label{equation:11}
\end{eqnarray}
where $M(s,t)$ is given in\eq{4}.
\end{proposition}
\begin{proof}
Recalling\eq{8}, we have
$$
\EE\left[X(t)\right]=1-\barF(t)\,\EE\left({\rm e}^{-W(t)}\right),
$$
from which Eq.\eq{10} immediately follows. Moreover, due to
\begin{eqnarray*}
&& \hspace{-1.7cm}
    \EE\left[X^2(t)\right]=\EE\left(1-2\barF(t)\,{\rm e}^{-W(t)}
    +\barF^2(t)\,{\rm e}^{-2 W(t)}\right) \\
&& \hspace{0.1cm}
    =1-2\barF(t)\EE\left({\rm e}^{-W(t)}\right)
    +\barF^2(t)\EE\left({\rm e}^{-2 W(t)}\right),
\end{eqnarray*}
we obtain:
\begin{eqnarray*}
&& \hspace{-1.7cm}
    \Var\left[X(t)\right]=1-2\barF(t)\EE\left({\rm e}^{-W(t)}\right)
    +\barF^2(t)\EE\left({\rm e}^{-2 W(t)}\right)
    \nonumber \\
&& \hspace{0.2cm}
    -1+2\barF(t)\EE\left({\rm e}^{-W(t)}\right)
    -\left[\EE\left({\rm e}^{-W(t)}\right)\right]^2\barF^2(t)
    \nonumber \\
&& \hspace{0.2cm}
    =\barF^2(t)\left\{\EE\left({\rm e}^{-2 W(t)}\right)
    -\left[\EE\left({\rm e}^{-W(t)}\right)\right]^2\right\},
\end{eqnarray*}
and thus Eq.\eq{11}.
\end{proof}
\begin{proposition}\label{prop:medvar}
The mean of $X(t)$ is increasing in $t\in\Stati$, and is decreasing in $c$. Moreover,
\begin{eqnarray}
&& \hspace{-1.7cm}
\lim_{t\to 0}\,\EE[X(t)]=0, \qquad \lim_{t\to 0}\,\Var[X(t)]=0;
\label{equation:14} \\
&& \hspace{-1.7cm}
\lim_{t\to \ell}\,\EE[X(t)]=1, \qquad  \lim_{t\to \ell}\,\Var[X(t)]=0.
\label{equation:15}
\end{eqnarray}
\end{proposition}
\begin{proof}
Using Eqs. (\ref{equation:4}) and (\ref{equation:20}),  
\begin{eqnarray*}
&& \hspace{-0.5cm}
 \frac{{\rm d}}{{\rm d}t}\left[\overline{F}(t) M(-1,t)\right] \\
&& \hspace{0.5cm}
=-\frac{r(t)+\lambda-\sqrt{\lambda^2+c^2}}{2}
\left(1+\frac{\lambda}{\sqrt{\lambda^2+c^2}} \right)
\exp\left\{-\int_0^t r(u) {\rm d}u -\lambda t
+\sqrt{\lambda^2 +c^2}\,t \right\}
\\
&& \hspace{0.5cm}
-\frac{r(t)+\lambda+\sqrt{\lambda^2+c^2}}{2}
\left(1-\frac{\lambda}{\sqrt{\lambda^2+c^2}} \right)
\exp\left\{-\int_0^t r(u) {\rm d}u -\lambda t
-\sqrt{\lambda^2 +c^2}\,t \right\}.
\end{eqnarray*}
Hence, recalling $r(t)>c$ for all $t\in{\cal S}$ we have 
${\frac{{\rm d}}{{\rm d}t}[\overline{F}(t) M(-1,t)]<0}$. 
By virtue of\eq{10} the mean of $X(t)$ is thus increasing in $t$. Moreover, 
from Eq.\eq{4} it is not hard to see that $M(s,t)$ is increasing in $c\geq 0$. 
Hence, $\EE[X(t)]$ is decreasing in $c$ due to\eq{10}. Limits\eq{14} easily follow 
from\eq{10} and\eq{11}. Recalling Eqs.\eq{10} and\eq{4} we have:
$$
\lim_{t\to \ell}\, \EE[X(t)]=1-\frac{1}{2} \left(1+\frac{\lambda}{\sqrt{\lambda^2+c^2}}\right)\,
\lim_{t\to \ell}\, \exp\left\{-\int_{0}^{t} r(u){\rm d}u-\lambda t+\sqrt{\lambda^2+c^2}\,t\right\},
$$
where
\begin{eqnarray}
&& \hspace{-1.2cm}
\lim_{t\to \ell}\, \exp\left\{-\int_{0}^{t} r(u){\rm d}u-\lambda t+\sqrt{\lambda^2+c^2}\,t\right\}
\nonumber
\\
&& \hspace{-0.2cm}
=\lim_{t\to \ell}\, \exp\left\{-\int_{0}^{t} \left[r(u)-c\right]\,{\rm d}u \right\}
\exp\left\{-\int_{0}^{t} \left[c+\lambda-\sqrt{\lambda^2+c^2}\right]\,{\rm d}u \right\}=0.
\label{equation:16}
\end{eqnarray}
The first of\eq{15} thus follows. Moreover, due to Eqs.\eq{11} and\eq{4}, we have:
\begin{eqnarray*}
 && \hspace{-0.4cm}
 \lim_{t\to \ell}\, \Var\left[X(t)\right] 
 =\frac{1}{2} \left(1+\frac{\lambda}{\sqrt{\lambda^2+4 c^2}}\right)\, \\
 && \hspace{0.8cm}
 \times \lim_{t\to \ell}\,\exp\left\{{-2\int_{0}^{t} [r(u)-c]\,{\rm d}u}\right\}
 \exp\left\{{-\int_{0}^{t} \left[2c+\lambda-\sqrt{\lambda^2+4 c^2}\right]\,{\rm d}u}\right\} \\
 && \hspace{0.8cm}
 -\frac{1}{4} \left(1+\frac{\lambda}{\sqrt{\lambda^2+c^2}}\right)^2 \\
 && \hspace{0.8cm}
 \times \lim_{t\to \ell}\,\exp\left\{{-2\int_{0}^{t} [r(u)-c]\,{\rm d}u}\right\}
 \exp\left\{{-2\int_{0}^{t} \left[c+\lambda-\sqrt{\lambda^2+c^2}\right]\,{\rm d}u}\right\}.
\end{eqnarray*}
%
%
Finally, noting that 
$$
\lim_{t\to \ell}\,\exp\left\{-2\int_{0}^{t} [r(u)-c]\,{\rm d}u\right\}
\exp\left\{-\int_{0}^{t} [2c+\lambda-\sqrt{\lambda^2+4 c^2}]\,{\rm d}u\right\}=0  
$$
and recalling Eq.\eq{16}, the second of\eq{15} holds.
\end{proof}
\par
Figure 4 shows mean and variance of $X(t)$ for an example involving distribution 
function\eq{22} and hazard rate\eq{23}.   
%
\begin{figure}[t] 
\begin{center}
\epsfxsize=8.5cm
\centerline{\epsfbox{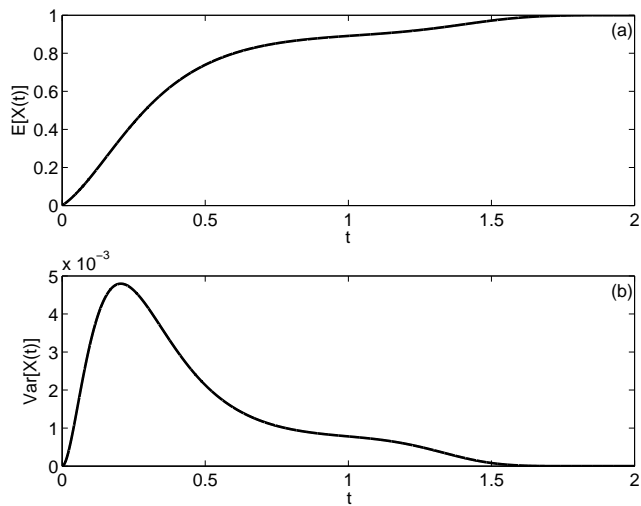}}
\caption{\small (a) Mean and (b) variance of $X(t)$; 
$r(t)$ and all parameters are as in case (a) of Figure 2, i.e.\ $c=1$, 
$\lambda=\alpha=15$, $\beta=0.001$ and $r(t)$ is given by Eq.\eq{23}.}
\label{fig:4}
\end{center}
\end{figure}
%
\par
An immediate consequence of Proposition \ref{prop:medvar} is the following
\begin{proposition}
Process $X(t)$ converges in probability to $1$ as $t\to \ell$. 
\end{proposition}
\par
Finally, denoting by $V(t)=V(0)\,(-1)^{N(t)}$ the telegraph process that 
describes the perturbating noise, making use of a well-known result (see, for instance, 
Theorem 3.4.4 of Ross \cite{Ro96}) the following asymptotic probabilities follow:
$$
 \lim_{t\to +\infty}\Prob\{V(t)=c\}=\lim_{t\to +\infty}\Prob\{V(t)=-c\}={1\over 2}.
$$
%
\section{\bf An application}\label{section:4}
In this section we indicate a statistical procedure based on an asymptotic confidence 
band for estimated hazard rates, which is useful to assess the validity of the 
proposed stochastic model. 
\par
Let $T_1,T_2,\ldots,T_n$ be iid absolutely continuous random variables 
describing a random sample of failure times having density $f(t)$ and
distribution function $F(t)$. The usual kernel estimator of $f(t)$ is 
\begin{equation}
 \widehat{f}(t)=\frac{1}{nh}\sum_{i=1}^{n} k\left(\frac{t-T_i}{h}\right),
 \label{equation:25}
\end{equation}
where $k(\cdot)$ is a bounded even density function, and $h$ is the bandwidth. 
The corresponding estimator of $F(t)$ is 
\begin{equation}
 \widehat{F}(t)=\int_{-\infty}^{t} \hat{f}(u) \,{\rm d}u
 =\frac{1}{n}\sum_{i=1}^{n} K\left(\frac{t-T_i}{h}\right),
 \label{equation:26}
\end{equation}
where $K(t)=\int_{-\infty}^t k(u) \,{\rm d}u$. From\eq{25} and\eq{26} we can 
build up the following estimator for the hazard rate function:
\begin{equation}
 \widehat{r}(t)=\frac{\widehat{f}(t)}{1-\widehat{F}(t)}.
 \label{equation:27}
\end{equation}
Taking into account the asymptotic normality of\eq{25} and\eq{26}, 
with $\Var(\widehat{f} )={(nh)^{-1}} \mathcal{K} f$, where 
${\mathcal{K}}=\int_{-\infty}^{+\infty} k^2(t)\,{\rm d}t$, and 
$\Var(\widehat{F})=O(n^{-1})$, it follows that also $\widehat{r}(t)$ is 
asymptotically normal, with variance $\Var(\widehat{r})={(nh)^{-1}} \mathcal{K} r^2/f$ 
(see Hall {\em et al.}\ \cite{Hall}). We can thus consider the following confidence 
bands for $r(t)$, having nominal coverage $1-2\alpha$, 
\begin{equation}
 \widehat{r}^{\,-}(t)=\widehat{r}(t)-\left[\ds\frac{\mathcal{K}}{nh\,\widehat{f}(t)}\right]^{1/2} 
 \widehat{r}(t)\,z_\alpha,
 \qquad 
\widehat{r}^{\,+}(t)=\widehat{r}(t)+\left[\ds\frac{\mathcal{K}}{nh\,\widehat{f}(t)}\right]^{1/2} 
 \widehat{r}(t)\,z_\alpha, 
 \label{equation:28}
\end{equation}
where $z_\alpha$ is the $\alpha$-level critical point of the standard normal distribution. 
\par
We note that a random hazard rate of type $r(t)+V(0)(-1)^{N(t)}$, related to process\eq{19}, 
alternates randomly between functions $r(t)-c$ and $r(t)+c$, thus being contained inside the 
strip $r(t)\pm c$. This suggests the following procedure to test if the random sample 
$(T_1,T_2,\ldots,T_n)$ can be viewed as drawn from the (random) distribution function $X(t)$: 
If, for a fixed $\alpha$, the strip $r(t)\pm c$ is contained within a realization of 
the confidence band\eq{28} obtained from observed data, the model defined by\eq{19} 
is defensible for the given data set. In other terms, with a $(1-2\alpha)$-level 
confidence, we can adopt model $X(t)$ if 
\begin{equation}
 |r(t)-\widehat{r}(t)|\leq 
 \left[\ds\frac{\mathcal{K}}{nh\,\widehat{f}(t)}\right]^{1/2}\widehat{r}(t)\,z_\alpha - c 
 \qquad \hbox{for all }t\in(t_{(1:n)},t_{(n-1:n)}), 
 \label{equation:29}
\end{equation}
where $t_{(j:n)}$ denotes the value of the $j$-th order statistic in the observed random sample. 
We stress that, due to\eq{29}, this procedure is effective to reject ``large'' values of $c$. 
%
\begin{figure}[t] 
\begin{center}
\epsfxsize=8.5cm
\centerline{\epsfbox{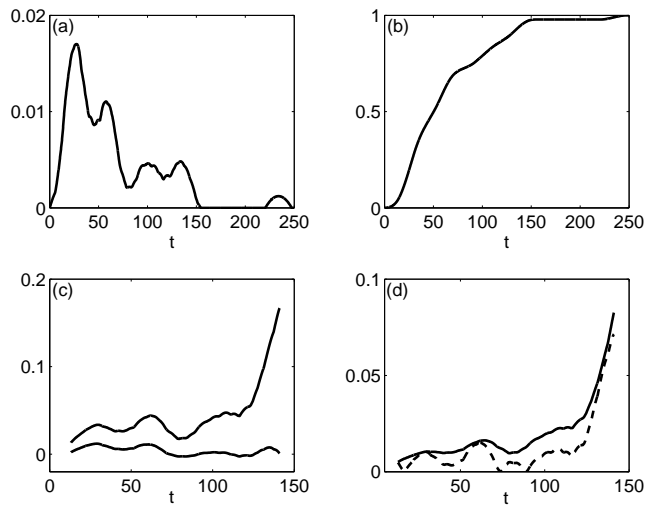}}
\caption{\small For the sample data of Application \ref{esempio1}, with $r(t)=0.0125$, $t\geq 0$, 
$\alpha=0.025$, $c=0.0004$, $h=6$ and $k(t)$ given in\eq{30}, it is shown:  
(a) the estimated density obtained by Eq.\eq{25}, 
(b) the estimated distribution function obtained by Eq.\eq{26}, 
(c) the realization of the confidence band\eq{28}, 
(d) the realization of the functions appearing in the left-hand-side (dashed line) 
and in the right-hand-side (solid line) of Eq.\eq{29}.}
\label{fig:5}
\end{center}
\end{figure}
%
\begin{application}\label{esempio1}{\rm
{\bf -- A case-study with constant hazard rate function.}\\ 
Let us consider the following set of $46$ sample data given in Ahmad \cite{Ah99}: 
\begin{eqnarray*}
&& \{13, 14, 19, 19, 20, 21, 23, 23, 25, 26, 26, 27, 27, 31, 32, 34, 34, \\
&& 37, 38, 38, 46, 46, 50, 53, 54, 57, 58, 59, 60, 65, 65, 66, 70, 85, \\
&& 90, 98, 102, 103, 110, 118, 124, 130, 136, 138, 141, 234\}. 
\end{eqnarray*}
These represent the survival times of certain patients in a melanoma study conducted 
by the Central Oncology Group (see Susarla and Van Ryzin \cite{SuVaRy86}). 
Figure \ref{fig:5}(a) shows the estimated density obtained by use of\eq{25}, 
where $k(t)$ is the Epanechnikov kernel (see Silverman \cite{Si86})
\begin{equation}
 k(t)=\cases{
   \ds\frac{3}{4 \sqrt{5}} \left(1-\ds\frac{t^2}{5}\right) & if $-\sqrt{5} \leq t \leq \sqrt{5}$\cr
   0 & otherwise.}
 \label{equation:30}
\end{equation}
The corresponding estimated distribution function is given in Figure \ref{fig:5}(b). 
A bandwidth depending on the sample standard deviation was proposed in Azzalini \cite{Az81}, 
where it is also mentioned that these results are hardly affected if the standard 
deviation is estimated from the sample data. Hence, we have empirically used the fixed 
bandwidth $h=6$; this choice has been also motivated by the need of obtaining a 
sufficiently smooth estimated density with small tails. We consider a baseline constant 
hazard rate function $r(t)=0.0125$, $t\geq 0$. The realization of the confidence 
band\eq{28} obtained from the above mentioned data is plotted in Figure \ref{fig:5}(c) 
for $\alpha=0.025$. With a $0.95$-level confidence the strip $r(t)\pm c$, with $c=0.0004$, 
falls within the confidence band, i.e.\ condition\eq{29} is fulfilled (this is graphically 
shown in Fig.\ \ref{fig:5}(d)). In conclusion, stochastic model\eq{19} is defensible 
for the observed data when $r(t)=0.0125$, $t\geq 0$, and $c\leq 0.0004$.  
}\end{application}
%
\begin{figure}[t] 
\begin{center}
\epsfxsize=8.5cm
\centerline{\epsfbox{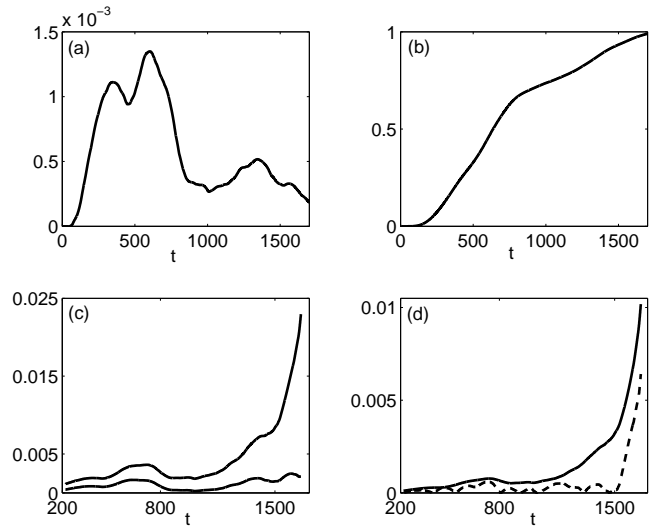}}
\caption{\small For the sample data of Application \ref{esempio2}, with $r(t)$ given in 
(\ref{equation:7}), $\alpha=0.025$, $c=0.00025$, $h=75$ and $k(t)$ given in\eq{30}, we show 
the same functions appearing in Figure \ref{fig:5}.}
\label{fig:6}
\end{center}
\end{figure}
%
\begin{application}\label{esempio2}{\rm
{\bf -- A case-study with varying hazard rate function.}\\ 
Hereafter we analyze the following set of $n=86$ sample data, taken from Table 1 of 
Langseth and Lindqvist \cite{LaLi05}: 
\begin{eqnarray*}
&& \{220, 233, 234, 240, 265, 270, 273, 279, 285, 287, 294, 295, 300, 325, 328,  \\
&& 333, 365, 368, 369, 381, 417, 418, 429, 460, 470, 474, 475, 476, 508, 522, 523,  \\
&& 535, 542, 570, 580, 604, 612, 613, 614, 615, 634, 636, 637, 638, 651, 657, 660,  \\
&& 666, 668, 680, 681, 684, 691, 693, 705, 717, 834, 837, 841, 843, 845, 875, 972,  \\
&& 1037, 1084, 1091, 1109, 1117, 1197, 1258, 1269, 1297, 1309, 1322, 1346, 1349,   \\
&& 1359, 1363, 1448, 1476, 1481, 1557, 1606, 1610, 1642, 1659\}. 
\end{eqnarray*}
They describes the service time of a single component in a reliability study. 
Figures \ref{fig:6}(a) and \ref{fig:6}(b) show the estimated density obtained from\eq{25}, 
with $k(t)$ given in\eq{30}, and the corresponding estimated distribution function. 
Here, the bandwidth appearing in\eq{25} has been empirically fixed as $h=75$, in agreement 
with the remarks in Application \ref{esempio1}. We choose the following baseline 
non-monotonic hazard rate function: 
\begin{equation}
 r(t)=\cases{
 3.5\cdot 10^{-6}\,t & if $0\leq t\leq 650$, \cr
 -4.07143\cdot 10^{-6}\,t+0.00492143 & if $650< t\leq 1000$, \cr
 8\cdot 10^{-6}\,t-0.00715 & if $t> 1000$.}
 \label{equation:7}
\end{equation}
Fig.\ \ref{fig:6}(c)  shows a realization of the confidence band\eq{28} obtained 
from the above data for $\alpha=0.025$. For $c=0.00025$ the strip $r(t)\pm c$ 
falls within the confidence band, with a $0.95$-level confidence. 
Condition\eq{29} is thus fulfilled (see Fig.\ \ref{fig:6}(d)). 
Hence, we are finally led to consider model\eq{19} as defensible 
for the observed data, with $r(t)$ given in\eq{7} and $c\leq 0.00025$.  
}\end{application}
%
\subsection*{\bf Acknowledgments}
This work has been performed under partial support by MIUR (PRIN 2005), 
by G.N.C.S.-INdAM and by Campania Region.

%
%
\end{document}